\documentclass{article}
\usepackage{amsmath}
\usepackage{amssymb}
\usepackage{graphicx}
\newtheorem{theorem}{Theorem}[section]
\newtheorem{Definition}[theorem]{Definition}
\newenvironment{definition}{\begin{Definition}\rm}{\end{Definition}}
\newtheorem{lemma}[theorem]{Lemma}
\newtheorem{proposition}[theorem]{Proposition}
\newtheorem{Remark}[theorem]{Remark}
\newenvironment{remark}{\begin{Remark}\rm}{\end{Remark}}
\newenvironment{proof}[1][Proof]{\noindent\textbf{#1.} }{\ \rule{0.5em}{0.5em}}

\numberwithin{equation}{section}

\newcommand\sLP{\\[\smallskipamount]}
\newcommand\mLP{\\[\medskipamount]}
\newcommand\bLP{\\[\bigskipamount]}
\newcommand\PP{\mathbb{P}}
\newcommand\RR{\mathbb{R}}
\newcommand\FSL{\mathcal{L}}
\newcommand\FST{\mathcal{T}}
\newcommand\al\alpha
\newcommand\be\beta
\newcommand\ga\gamma
\newcommand\de\delta
\newcommand\om\omega
\newcommand\La\Lambda
\newcommand\Om\Omega
\newcommand\iy\infty
\newcommand\half{\frac12}
\newcommand\thalf{\tfrac12}
\newcommand{\hyp}[5]{\,\mbox{}_{#1}F_{#2}\!\left(
  \genfrac{}{}{0pt}{}{#3}{#4};#5\right)}
\newcommand\Lcrit{\FSL_{\operatorname{critical}}}

\begin{document}

\title{The Intersection of Bivariate Orthogonal Polynomials on Triangle Patches}
\author{Tom H. Koornwinder\thanks{Korteweg-de Vries Institute for Mathematics,
University of Amsterdam, P.O.\ Box 94248, 1090 GE Amsterdam, Netherlands,
e-mail: \texttt{T.H.Koornwinder@uva.nl }}
\and Stefan A. Sauter\thanks{Institut f\"{u}r Mathematik, Universit\"{a}t
Z\"{u}rich, Winterthurerstrasse 190, CH-8057 Z\"{u}rich, Switzerland, e-mail:
\texttt{stas@math.uzh.ch}}}
\date{}
\maketitle

\begin{abstract}\small
In this paper, the intersection of bivariate orthogonal polynomials
on triangle patches will be investigated. The result is interesting
by its own but also has important applications in the theory
of a posteriori error estimation for finite element discretizations
with $p$-refinement, i.e., if the local polynomial degree of
the test and trial functions is increased to improve the accuracy.
A triangle patch is a set of disjoint open triangles whose closed union covers
a neighborhood of the common triangle vertex.
On each triangle we consider the space of orthogonal polynomials of degree $n$
with respect to the weight function which is the product of the barycentric
coordinates.
We show that
the intersection of these polynomial spaces is the null space.
The analysis requires the derivation of subtle
representations of orthogonal polynomials on triangles.
Up to four triangles have to be considered to identify
that the intersection is trivial.
\mLP
\textbf{Keywords:}
A posteriori error estimation,
saturation property,
$p$-refinement,
Jacobi polynomials,
orthogonal polynomials on the triangle,
intersections of $n$-th degree orthogonal polynomial spaces
\mLP
\textbf{Mathematics Subject Classification (2010):}
65N15, 65N30, 65N50, 33C45, 33C50
\end{abstract}

\section{Introduction}
In this paper, we will investigate the intersection of bivariate orthogonal
polynomials on triangle patches. This problem arises in the theory of a
posteriori error estimation for finite element discretizations of elliptic
partial differential equations --- in particular if the local polynomial degree
of the finite element spaces is increased during the solution process. Before
we give the precise mathematical formulation of this problem we will sketch
its application in the finite element analysis.

A posteriori error estimation and adaptivity are well established
methodologies for the numerical solution of partial differential equations by
finite elements (cf. \cite{BaRe1}, \cite{BaRe2}, \cite{Verfuerth},
\cite{AiOd}, \cite{BangRa}, \cite{Repin},
\cite{Doerfler}, \cite{MekNoch}, \cite{Stevenson},
\cite{BiDaDeV}).

Some types of error estimators, for example hierarchical ones (see,
e.g., \cite{BankWei}, \cite{BornErdKorn}, \cite{BankXu}), require
explicitly or implicitly the \emph{saturation} assumption which states that
the error on the refined mesh and/or with higher polynomial degree is strictly
smaller than the error on the previous mesh/polynomial degree. In the
pioneering paper \cite{DoerflerNoch} the saturation assumption is proved
for the Poisson problem in two spatial dimensions under the assumption that
the data oscillations are small. In \cite{MekNoch} the convergence of
adaptive finite element methods (AFEM) for general (nonsymmetric) second order
linear elliptic partial differential equations is proved, where the term
``adaptivity'' is understood in the sense of
adaptive \emph{mesh} refinement and the polynomial degree stays fixed. The
theory in \cite{MekNoch} also generalizes the proof of the saturation
property to quite general 2nd order elliptic problems and estimates the error
on the refined mesh by the error of the coarser mesh plus a data oscillation term.

For the proof of the saturation assumption for $p$-refinement, i.e., when the
local polynomial degree of the finite element space is increased instead of
the mesh being refined, a difficulty arises which is related to a polynomial
projection property on triangle patches. Here orthogonal polynomials in
two variables on a triangle (see \cite{Proriol}, 
\cite{Koo},  \cite{DunklXu})
enter, and the problem just raised
is also interesting by
its own in that area.
By the way,
these orthogonal polynomials also have important applications in the
field of spectral methods
for discretizing partial differential equations and we refer to \cite{HestWar}
for further details. In
particular, orthogonal polynomials on triangles can be efficiently used for
discontinous Galerkin (dG) methods
or to discretize boundary integral equation of negative order since no
continuity is required across simplex boundaries.

\begin{figure}[ptb]
\begin{center}
\includegraphics{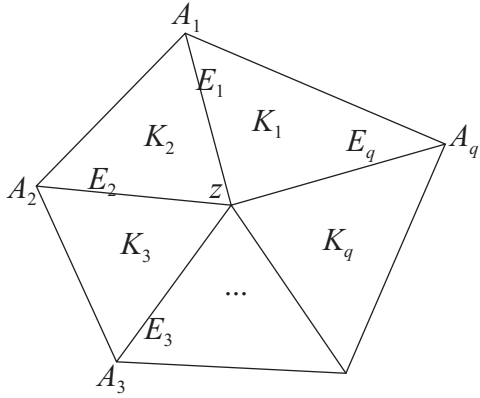}
\caption{Triangle patch $\FST :=\left\{  K_{i}:1\leq i\leq q\right\}  $
around a point $z$. The triangles $K_{i}$ and $K_{i+1}$ share an edge $E_{i}$.
Each triangle $K_{i}$ has $z$ as a
vertex while its other vertices are denoted by $A_{i-1}$ and $A_{i}$.}
\label{Fig1}
\end{center}
\end{figure}

Let us now briefly state the problem which we will solve in this paper.
For a two-dimensional domain $D\subset\mathbb{R}^{2}$, the set of bivariate
polynomials of maximal total degree $n$ is denoted by
$\PP_{n}\left(D\right)$. Put $\PP_{-1}(D):=\{0\}$.
Let $z\in\RR^{2}$ and let $\FST :=
\left\{K_{i}:1\leq i\leq q\right\}  $ denote a \emph{triangle patch around }$z$,
i.e., $\FST $ is a set of (open) triangles (cf.~Figure \ref{Fig1}) which
\begin{itemize}
\item are pairwise disjoint,
\item have $z$ as a common vertex.
\item For all $1\leq i\leq q$, the triangles\footnote{We use here the
convention $K_{0}:=K_{q}$ and $K_{q+1}:=K_{1}$ and analogously for
the vertices $A_{i}$ and the edges
$E_{i}$. Clearly $q\geq3$ holds.} $K_{i}$ and $K_{i+1}$
have one common edge, denoted by $E_{i}$, which connects the common
vertices $z$ and $A_i$ of $K_i$ and $K_{i+1}$.
\end{itemize}
Thus
\begin{equation}
K_{i}=\operatorname*{conv^o}(z,A_{i-1},A_{i}),
\label{K15}
\end{equation}
where conv denotes the convex hull of the given points and
$\operatorname*{conv^o}$ the open interior of this convex hull.

Let
$\Om:=\operatorname*{int}
\left(\cup_{i=1}^{q}\overline{K_{i}}\,\right)$
and let
$\mathfrak{S}:=\Om\cap\left(\cup_{i=1}^{q}\partial K_{i}\right)$
denote the inner mesh skeleton. We denote by
$\PP_{n}\left(\FST \right)$
the space of piecewise polynomials, i.e.,%
\[
\PP_{n}\left(  \FST \right)  :=\left\{  f:\Om\backslash
\mathfrak{S}\rightarrow\RR\mid\forall i\in\{1,\ldots,q\}\;\;
f\big|_{K_{i}}\in\PP_{n}\left(  K_{i}\right)  \right\}  .
\]
We consider $\PP_n(\Om)$ as a linear subspace of $\PP_n(\FST)$ by its
natural embedding.

For $i=1,\ldots,q$ define a weight function $\om_i:=a_ib_ic_i$ on $K_i$,
where $a_i,b_i,c_i$ are affine linear functions
which vanish on the respective edges of $K_i$.
Thus $\om_i$ is the product of the barycentric coordinates in $K_{i}$ or,
in other terms, a cubic bubble function which is positive on $K_{i}$.
We define the inner product $\left(\cdot,\cdot\right)_{\FST }$
on $\FST $ by
\begin{equation}
\left(u,v\right)_{\FST}:=\sum_{i=1}^{q} (u,v)_{K_i},
\label{K25}
\end{equation}
where
\begin{equation}
(u,v)_{K_i}:=\int_{K_{i}}
u(x,y)\,v(x,y)\,\om_i(x,y)\,dx\,dy.
\label{K31}
\end{equation}
Denote by $\PP_{n-1}^\perp(K_i)$ the orthoplement of
$\PP_{n-1}(K_i)$ in $\PP_n(K_i)$ with respect to the inner product
\eqref{K31}.
Let $\Pi_n^{\FST}\colon\PP_{n}(\Om)\to\PP_{n-1}(\FST)$
denote the restriction to $\PP_n(\Om)$ of the orthogonal projection
of $\PP_n(\FST)$ onto $\PP_{n-1}(\FST)$ with respect to the inner product \eqref{K25}.

\begin{theorem}
\label{MainConj}Let $n\geq1$.
Then the following three statements are equivalent and each of them holds.
\mLP
{\bf(a)} If $u\in\PP_n(\Om)$ and $(u,w)_\FST=0$ for all $w\in\PP_{n-1}(\FST)$
then $u=0$.
\mLP
{\bf(b)}
$\displaystyle\bigcap_{i=1,\ldots,q} \PP_{n-1}^\perp(K_i)=\{0\}$.
\mLP
{\bf(c)} The map $\Pi_n^{\FST}\colon\PP_{n}(\Om)\to\PP_{n-1}(\FST)$
 is injective.
\end{theorem}

The equivalence of the three statements is trivial, so we can pick one of them as what we aim to prove.
It turns out that (b) is the most convenient statement for a proof.
Then it is natural to examine first the intersection of two such spaces
for adjacent triangles, i.e.,
$\PP_{n-1}^{\perp}(K_i)\cap\PP_{n-1}^{\perp}(K_{i-1})$.
This will be the subject of Section \ref{K32}, where
explicit knowledge of orthogonal polynomials on the triangle for the inner product
\eqref{K31}, to be summarized in Section \ref {sec_OPtriangle}, is crucial.

By Section \ref{K32} the intersection for two adjacent triangles is mostly $\{0\}$,
but there
are exceptional cases. For these cases it is necessary to consider the intersection
of spaces for three adjacent triangles, and in one case for four adjacent triangles,
in order to get an intersection $\{0\}$. This is the subject of
Section \ref{Secn>1} (for $n>1$) and of Section \ref{Secn=1} (for $n=1$).

The equivalent formulation (c) of Theorem \ref{MainConj} raises the question
to estimate $(v,\Pi_n^\FST v)_\FST$ from below, also in dependence of the
triangle patch $\FST$. Some generalities about this  will be given in the
final Section \ref{K33}.

In principle, all computations in this paper can be done by hand.
Nevertheless, some of the more tedious computations we have done in
\textit{Mathematica}${}^{\mbox{\footnotesize\textregistered}}$,
while we have also checked many of the other computations by this program.

It is quite probable that the results and proofs in this paper can be carried over to the case that $\om_i:=(a_ib_ic_i)^\al$ ($\al>-1$), i.e., that the weight function
is some power of the product of the barycentric coordinates. We have refrained
from doing the computations in this more general case because only the special
case is needed in the application we have sketched.

\section{Orthogonal polynomials on the triangle}
\label{sec_OPtriangle}
Let $\al,\be>-1$.  The \emph{Jacobi polynomial} $P_n^{(\al,\be)}$ 
(see for instance \cite{Szego}) is a polynomial of degree $n$ such that
\[
\int_{-1}^1 P_n^{(\al,\be)}(x)\,q(x)\,(1-x)^\al(1+x)^\be\,dx=0\
\]
for all polynomials $q$ of degree less than $n$, and
\[
P_n^{(\al,\be)}(1)=\frac{(\al+1)_n}{n!}\,.
\]
Here the \emph{shifted factorial} is defined by $(a)_n:=a(a+1)\ldots(a+n-1)$ for $n>0$
and $(a)_0:=1$.
All zeros of $P_n^{(\al,\be)}$ lie in $(-1,1)$, so it has definite sign
on $[1,\infty)$ and on $(-\infty,-1]$.

The Jacobi polynomial has an explicit expression in terms of
a \emph{terminating Gauss hypergeometric series}
\[
\hyp21{-n,b}cz:=\sum_{k=0}^n\frac{(-n)_k (b)_k}{(c)_k\,k!}\,z^k
\]
as follows:
\[
P_n^{(\al,\be)}(x)=\frac{(\al+1)_n}{n!}\,\hyp21{-n,n+\al+\be+1}{\al+1}{\frac{1-x}2}.
\]
There is the symmetry relation
\[
P_n^{(\al,\be)}(-x)=(-1)^nP_n^{(\be,\al)}(x).
\]

Let $T_1$ be the open (unit) triangle
\begin{equation}
T_1:=\{(x,y)\mid x,y,1-x-y>0\}.
\label{K12}
\end{equation}
Let $\al,\be,\ga>-1$.
Define in terms of Jacobi polynomials the bivariate polynomial
\begin{equation}
P_{n,k}^{(\al,\be,\ga)}(x,y):=
(1-x)^k\,P_{n-k}^{(\al,\be+\ga+2k+1)}(1-2x)\,
P_k^{(\be,\ga)}\left(1-\frac{2y}{1-x}\right).
\label{K2}
\end{equation}
This is a polynomial of degree $n$ in $x$ and $y$.
For $(n,k)\ne(m,j)$ we have the orthogonality relation
\begin{equation}
\int_{T_1} P_{n,k}^{(\al,\be,\ga)}(x,y) P_{m,j}^{(\al,\be,\ga)}(x,y)\,
w_{\al,\be,\ga}(x,y)\,dx\,dy=0,
\label{K1}
\end{equation}
where $w_{\al,\be,\ga}(x,y):=x^\al y^\be (1-x-y)^\ga$.
This follows immediately from the orthogonality relations for Jacobi polynomials
if we write
\begin{multline*}
\int_{T_1} f(x,y)\,dx\,dy=\int_0^1\left(\int_0^{1-x} f(x,y)\,dy\right)dx\\
=\int_0^1(1-x)\left(\int_0^1 f(x,(1-x)t)\,dt\right)dx.
\end{multline*}
Thus with respect to
the inner product for $L^2(T_1,x^\al y^\be(1-x-y)^\ga\,dx\,dy)$
the system
$\{P_{m,j}^{(\al,\be,\ga)}\}_{0\le j\le m\le n}$
is an orthogonal basis of $\PP_n(T_1)$.

These bivariate orthogonal polynomials on the triangle were introduced
by Proriol \cite{Proriol}, see also the survey
\cite{Koo} and the monograph \cite{DunklXu}.
In the context of numerical analysis they were rediscovered in special
cases in \cite{Dubiner} and they got ample coverage in the monograph
\cite{KarSher}.

Denote by $\PP_{n-1}^\perp(T_1)$ the orthoplement of
$\PP_{n-1}(T_1)$ in $\PP_n(T_1)$ with respect to the inner product
just mentioned.
(So $\PP_{-1}^\perp(T_1)=\PP_0(T_1)$ consists of the constant functions.)$\;$
Then the system $\{P_{n,j}^{(\al,\be,\ga)}\}_{0\le j \le n}$ is a basis of
$\PP_{n-1}^\perp(T_1)$.
In particular, the polynomial
$(x,y)\mapsto P_n^{(\al,\be+\ga+1)}(1-2x)$
is in $\PP_{n-1}^\perp(T_1)$.

The symmetric group $S_3$ naturally acts on $T_1$. By considering the action
of $S_3$
on \eqref{K1} we obtain five further orthogonal bases for $\PP_n$ with
respect to the inner product for $L^2(T_1,x^\al y^\be(1-x-y)^\ga\,dx\,dy)$.
The six bases are as follows (considered as functions
of $(x,y)$).
\begin{equation}
\label{K3}
\begin{split}
&\{P_{m,j}^{(\al,\be,\ga)}(x,y)\}_{0\le j\le m\le n}\,,\qquad\qquad\;\,
\{P_{m,j}^{(\be,\al,\ga)}(y,x)\}_{0\le j\le m\le n}\,,\\
&\{P_{m,j}^{(\be,\ga,\al)}(y,1-x-y)\}_{0\le j\le m\le n}\,,\quad
\{P_{m,j}^{(\al,\ga,\be)}(x,1-x-y)\}_{0\le j\le m\le n}\,,\\
&\{P_{m,j}^{(\ga,\al,\be)}(1-x-y,x)\}_{0\le j\le m\le n}\,,\quad
\{P_{m,j}^{(\ga,\be,\al)}(1-x-y,y)\}_{0\le j\le m\le n}\,.
\end{split}
\end{equation}
In particular, each of these systems, when only taken for $m=n$, $0\le j\le n$,
is an orthogonal basis for $\PP_{n-1}^\perp(T_1)$.
In combination with \eqref{K2} this shows that the following polynomials in $(x,y)$
are elements of $\PP_{n-1}^\perp(T_1)$:
\begin{equation}
P_n^{(\al,\be+\ga+1)}(1-2x),\quad
P_n^{(\be,\al+\ga+1)}(1-2y),\quad
P_n^{(\ga,\al+\be+1)}(2(x+y)-1).
\label{K11}
\end{equation}

If $T$ is another open triangle in $\RR^2$ and
if $\La$ is an affine transformation of $\RR^2$ which maps $T$ onto $T_1$ then
the polynomials $P_{n,k}^{(\al,\be,\ga)}\circ\La$ are orthogonal on $T$
with respect to the weight function $w_{\al,\be,\ga}\circ\La$.
If $\al=\be=\ga$ then the inner product on $T$ is independent, up to constant factor,
of the choice of $\La$. In the sequel we will have $\al=\be=\ga=1$.
Similarly as for $T_1$, we denote by $\PP_{n-1}^\perp(T)$ the orthoplement
of $\PP_{n-1}(T)$ in $\PP_n(T)$ with respect to this inner product.

\section{The intersection of $n$-th degree orthogonal polynomial spaces
for two adjacent triangles}
\label{K32}
In this section we keep using the conventions and definitions of 
Section \ref{sec_OPtriangle}
for $\al=\be=\ga=1$, and we will compare the orthogonal polynomials
on the triangle $T_1$ for the weight function $w_{1,1,1}$ with
the orthogonal polynomials on the adjacent triangle
\begin{equation}
K_{c,d}:=\operatorname*{conv^o}\left( (1,0),(0,0),
\left(\tfrac{-c}{d-c},\tfrac1{d-c}\right)\right)\qquad(c\ne d)
\label{defKparab}
\end{equation}
for the weight function $w_{1,1,1}\circ\La$, where $\La$ is the affine map
sending $K_{c,d}$ to $T_1$, which is given by
\[
\La(x,y)=(x+cy,(d-c)y).
\]
We will prove:
\begin{theorem}
\label{Theo2Triangles}
For $n>2$ the intersection of the spaces of orthogonal
polynomials of
degree $n$ on $T_1$ and $K_{c,d}$, i.e., the space
$\PP_{n-1}^\perp(T_1)\cap \PP_{n-1}^\perp(K_{c,d})$,
has dimension zero unless $c=0$ or $d=1$
or $d-c=1$ or $c=1$, $d=0$. If $c=0$, $d=1$ then the intersection trivially
has dimension $n+1$. In all other exceptional cases the intersection has
dimension $1$.

For $n=2$ the intersection has dimension zero unless $c=0$ or $d=1$
or $d-c=\pm1$. If $c=0$, $d=1$ then the intersection trivially
has dimension $3$. In all other exceptional cases the intersection has
dimension $1$.

For $n=1$ the intersection has dimension $1$ except in the trivial case
$c=0$, $d=1$, when it has dimension $2$.

In the cases that the intersection has dimension $1$, it is spanned by
a polynomial $q_n^{(c,d)}(x,y)$ as follows:
\begin{align}
q_n^{(0,d)}(x,y)&=P_n^{(1,3)}(1-2x),\label{K18}\\
q_n^{(c,1)}(x,y)&=P_n^{(1,3)}(2(x+y)-1),\label{K19}\\
q_n^{(c,c+1)}(x,y)&=P_n^{(1,3)}(1-2y),\label{K20}\\
q_n^{(1,0)}(x,y)&=y^{-1}\left(  P_{n+1}^{(1,1)}(1-2(x+y))-P_{n+1}^{(1,1)}(1-2x)\right),\
\label{K16}\\
q_2^{(c,c-1)}(x,y)&=28\big(6x^2+6cxy+c(c+1)y^2\big)\nonumber\\
&\qquad-21(c+3)(2x+cy)+3c^2+15c+18,\label{K21}\\
q_1^{(c,d)}(x,y)&=3(c-d+1)x+3cy-2c+d-1.\label{K22}
\end{align}
\end{theorem}

In the case $d-c=1$ the triangles $T_{1}$ and $K_{c,d}$ have nonempty open
intersection, so for the application we have in mind the result for this case
is not needed.

Observe that, for $n=2$, \eqref{K21} agrees up to a constant factor
with \eqref{K18}, \eqref{K19}, \eqref{K16} for $c=0,2,1$, respectively.

For usage in the proof we pick from the orthogonal systems in \eqref{K3} one
particular orthogonal basis for $\PP_{n-1}^\perp(T_1)$, and we also renormalize it.
The resulting basis consists of the following polynomials $p_{n,k}$ ($k=0,\ldots,n$).
\begin{align}
p_{n,k}(x,y)&:=\frac{P_{n,k}^{(1,1,1)}(y,x)}{P_{n,k}^{(1,1,1)}(0,0)}\label{K4}
\\
&\;=\mbox{}_{2}F_{1}\!\left(\genfrac{}{}{0pt}{}{-n+k,n+k+5}{2};y\right)  (1-y)^{k}\,
\mbox{}_{2}F_{1}\!\left(\genfrac{}{}{0pt}{}{-k,k+3}{2};\frac{x}{1-y}\right)\nonumber\\
&\;=\mbox{}_{2}F_{1}\!\left(\genfrac{}{}{0pt}{}{-n+k,n+k+5}{2};y\right)
\sum_{j=0}^{k}\frac{(-k)_{j}(k+3)_{j}}{(2)_{j}j!}x^{j}(1-y)^{k-j}.\label{K10}
\end{align}
Similarly, for an orthogonal basis of $\PP_{n-1}^\perp(K_{c,d})$ we will take
the polynomials $(x,y)\mapsto p_{n,k}(x+cy,(d-c)y)$ ($k=0,\ldots,n$).
\\[\bigskipamount]\noindent
\textbf{Proof of Theorem \ref{Theo2Triangles}. }
By \eqref{K11} the polynomials
\begin{equation}
P_n^{(1,3)}(1-2x),\quad
P_n^{(1,3)}(1-2y),\quad
P_n^{(1,3)}(2(x+y)-1)
\label{K13}
\end{equation}
in $(x,y)$ are elements of $\PP_{n-1}^\perp(T_1)$, and the polynomials
\begin{equation}
P_n^{(1,3)}(1-2(x+cy)),\quad
P_n^{(1,3)}(1-2(d-c)y),\quad
P_n^{(1,3)}(2(x+dy)-1)
\label{K14}
\end{equation}
are elements of $\PP_{n-1}^\perp(K_{c,d})$. Hence for $c=0$ or 
$d=1$ or $d-c=1$ the intersection considered in the Theorem has dimension
at least 1, and the polynomial given by \eqref{K18}, \eqref{K19}, \eqref{K20},
respectively, is in this intersection.

A general element in $\PP_{n-1}^{\perp}\left(T_1\right)$ has the
form $\sum_{k=0}^n\alpha_k p_{n,k}(x,y)$, and a general element in
$\PP_{n-1}^{\perp}\left(K_{c,d}\right)$ has the form
$\sum_{k=0}^n\beta_k p_{n,k}(x+cy,(d-c)y)$.
Hence each nonzero element in $\PP_{n-1}^{\perp}\left(T_1\right)
\cap\PP_{n-1}^{\perp}\left(K_{c,d}\right)$
corresponds to a nontrivial solution of the
homogeneous linear system of equations
\begin{equation}
\begin{split}
&\mathrm{coeff}\left(  \sum_{k=0}^{n}\big(\alpha_{k}p_{n,k}(x,y)-\beta
_{k}p_{n,k}(x+cy,(d-c)y)\big),x^{r}y^{m}\right)=0\\
&\qquad\qquad\qquad\qquad\qquad\qquad\qquad\qquad
(r,m=0,1,\ldots,n,\;r+m\leq n)
\end{split}
\label{K5}
\end{equation}
of $\thalf(n+1)(n+2)$ equations in the $2(n+1)$ unknowns $\alpha_{0},\ldots
,\alpha_{n},\beta_{0},\ldots,\beta_{n}$.

By \eqref{K4} and \eqref{K2} the $n+1$ equations in \eqref{K5}
involving the coefficient of $x^r$ ($r=0,\ldots,n$) amount to
\[
\sum_{k=0}^n (\al_k-\be_k)\,
\frac{P_k^{(1,1)}(1-2x)}{P_k^{(1,1)}(1)}=0,
\]
which implies $\al_k=\be_k$ ($k=0,\ldots,n$).
So the system of equations \eqref{K5} reduces to
\begin{equation}
\begin{split}
&\sum_{k=0}^{n}\alpha_k\,\mathrm{coeff}\left( p_{n,k}(x,y)
-p_{n,k}(x+cy,(d-c)y),x^{r}y^{m}\right)  =0\\
&\qquad\qquad\qquad\qquad\qquad\qquad
(r=0,\ldots,n-1,\;m=1,\ldots,n,\;r+m\leq n),
\end{split}
\label{K7}
\end{equation}
which are $\thalf  n(n+1)$ homogeneous linear equations in the $n+1$
unknowns $\al_0,\ldots,\al_n$.

First we consider the case $n=1$. From \eqref{K10} we get that
\[
p_{1,0}(x,y)=1-3y\quad
p_{1,1}(x,y)=1-2x-y.
\]
Then we have to solve $\al_0,\al_1$ from the equation
\[
\al_0(1-3y)+\al_1(1-2x-y)=
\al_0(1-3(d-c)y)+\al_1(1-2x-(d+c)y).
\]
This yields
\[
3(d-c-1)\al_0+(d+c-1)\al_1=0,
\]
which has (if not $c=0$, $d=1$) a one-dimensional solution space spanned
by $(\al_0,\al_1):=(\thalf(1-d-c),\tfrac32(d-c-1))$. Then
$q_1^{(c,d)}(x,y)$ given by \eqref{K22} equals
$\al_0 p_{1,0}(x,y)+\al_1 p_{1,1}(x,y)$.

Now let $n\ge2$.
The power series coefficients in the left-hand sides of the equations \eqref{K7}
can be computed by using \eqref{K10}. We can rewrite the system \eqref{K7} as
\begin{equation}
\frac1{r!}\sum_{k=r}^n \al_k f_{k,r,m}(c,d)=0\qquad
(r=0,\ldots,n-1,\;m=1,\ldots,n,\;r+m\leq n),
\label{K8}
\end{equation}
where
\begin{align*}
f_{k,r,m}(c,d)&=
\sum_{i=\max(0,m-k+r)}%
^{\min(m,n-k)}\frac{(-n+k)_{i}(n+k+5)_{i}}{(2)_{i}\,i!}\\
&\quad\qquad\times
\frac{(-k)_{r}(k+3)_{r}(r-k)_{m-i}}{(2)_{r}(m-i)!}\,(1-(d-c)^{m})\\
&\quad-\sum_{i=\max(0,m-k+r)}^{\min(m-1,n-k)}
\frac{(-n+k)_{i}(n+k+5)_{i}}{(2)_{i}\,i!}\\
&\quad\qquad\times\sum_{j=r+1}^{\min(k,m+r-i)}
\frac{(-k)_{j}(k+3)_{j}(j-k)_{m+r-i-j}}{(2)_{j}(j-r)!(m+r-i-j)!}\,
c^{j-r}(d-c)^{m+r-j}.
\end{align*}
In particular,
\[
f_{r,r,m}(c,d)=
\frac{(-1)^{r}(2r+2)!}{(r+1)\,(r+2)!}\,\frac{(-n+r)_{m}\,(n+r+5)_{m}}
{(2)_{m}\,m!}\,(1-(d-c)^{m}),
\]
which is nonzero (note that $0\le m\le n-r$, hence $(-n+r)_m\ne0$)
except if $d=c+1$ or $d=c-1$ and $m$ even. Thus, for $d\ne c+1$
the system \eqref{K8} has a subsystem
\begin{equation}
\sum_{k=r}^n \al_k f_{k,r,1}(c,d)=0\qquad(r=0,\ldots,n-1)
\label{K9}
\end{equation}
with $f_{r,r,1}(c,d)\ne0$. 
Thus $\al_n$ successively determines $\al_{n-1},\al_{n-2},\ldots,\al_0$, by which
the system \eqref{K8} has a solution space of dimension at most 1 if
$d\ne c+1$.

For $d=c+1$ we have $f_{r,r,m}(c,c+1)=0$, while
\begin{equation*}
f_{r+1,r,m}(c,c+1)=
(-1)^{m+r+1}\,\frac{(n-r-1)!\,(n+r+6)_{m-1}}{(n-r-m)!\,(m-1)!\,m!}\,
\frac{(r+4)_{r+1}}{r+2}\,c.
\end{equation*}
This is nonzero unless $c=0$, but if $c=0$ then $d=1$ and we are in the
trivial case. Thus, for $d=c+1$, $c\ne0$ we successively get from \eqref{K9} together with
$f_{r,r,m}(c,c+1)=0$, $ f_{r+1,r,m}(c,c+1)\ne0$
that $\al_n=0$, $\al_{n-1}=0$,\ldots,$\al_1=0$. So only $\al_0$ may be nonzero
by which the system \eqref{K8} has a solution space of dimension at most 1.
In the beginning of the Proof we already saw that this dimension is at least 1.
This settles the case $d=c+1$ in the Theorem.

In the next step we consider the cases
$(r,m)=(n-1,1),(n-2,1),(n-2,2)$ of \eqref{K8}
(for $n=2$ these are all possible cases). This gives the following
system of three homogeneous linear equations in $\al_n,\al_{n-1},\al_{n-2}$.
\begin{equation}
\begin{split}
&  -\,\frac{(-1)^{n}(c+d-1)(2n+1)!}{(n-1)!\,(n+2)!}\,\alpha_{n}+\frac
{2(-1)^{n}(c-d+1)(n+2)(2n-1)!}{(n-1)!\,(n+1)!}\,\alpha_{n-1}=0,\\
&  \frac{(-1)^{n}(c+(c+d-1)n)(2n)!}{(n-2)!\,(n+2)!}\,\alpha_{n}\\
&  \qquad-\,\frac{(-1)^{n}\,2^{2n-1}(c(n+1)-(d-1)(n+3))(\tfrac{1}{2})_{n}%
}{n(n+1)(n-2)!}\,\alpha_{n-1}\\
&  \qquad-\,\frac{2(-1)^{n}(c-d+1)(2n+3)(2n-3)!}{(n-2)!\,n!}\,\alpha
_{n-2}=0,\\
&  -\,\frac{(-1)^{n}\big(2cd+((c+d)^{2}-1)n\big)(2n)!}{2(n-2)!\,(n+2)!}%
\,\alpha_{n}\\
&  \qquad+\frac{(-1)^{n}(c^{2}-d^{2}+1)(n+2)(2n-1)!}{(n-2)!\,(n+1)!}%
\,\alpha_{n-1}\\
&  \qquad-\,\frac{2(-1)^{n}((c-d)^{2}-1)(n+2)(2n+3)(2n-3)!}{3(n-2)!\,n!}%
\,\alpha_{n-2}=0.
\end{split}
\label{K6}
\end{equation}
The $3\times 3$ determinant of the coefficients of the system \eqref{K6}
can be computed to be equal to
\[
\frac{(-1)^{n+1}2^{4n}(2n)!\,(\tfrac{1}{2})_{n-1}(\tfrac{1}{2})_{n+2}%
}{3(n-2)!\,n!\,((n+1)!)^{2}}\,c(d-1)(c-d+1)(c-d-1).
\]
Thus $\alpha_{n}=\alpha_{n-1}=\alpha_{n-2}=0$ if not $c=0$ or $d=1$ or
$c-d=\pm1$.
Together with \eqref{K9} and $f_{r,r,1}(c,d)\ne0$ this implies that all $\al_k$ are
zero if not $c=0$ or $d=1$ or $c-d=\pm1$. This settles the Theorem in
the non-exceptional case except if $c-d=1$.
For $c=0$ or $d=1$ the Theorem is also settled now
because we already observed in the beginning of the Proof that the solution space
has dimension at least 1 in these cases.

Now consider the case $d=c-1$. Then the third equation in \eqref{K6}
is a multiple of the first equation, so we can solve from the first and second
equation of \eqref{K6} that
\begin{equation}
\begin{split}
&  \alpha_{n-1}=\frac{(c-1)n(2n+1)}{(n+2)^{2}}\,\alpha_{n},\\
&  \alpha_{n-2}=\frac{(n-1)n(2n-1)\big(3(n+1)+c(c-2)(2n+1)\big)}
{(n+1)(n+2)^{2}(2n+3)}\,\alpha_{n}.
\end{split}
\label{K17}
\end{equation}
For $n=2$ we conclude that the intersection has dimension 1
and that it contains the polynomial
$p_{2,2}(x,y)+\alpha_1 p_{2,1}(x,y)+\alpha_0 p_{2,0}(x,y)$ with $\alpha_1$
and $\alpha_0$ given by \eqref{K17} for $n=2$ and $\alpha_2=1$.
Together with \eqref{K10} this yields \eqref{K21}.

For $n>2$ we plug the above two equations into the cases
$(r,m)=(n-3,1),(n-3,2),(n-3,3)$ of \eqref{K8} with $d=c-1$. There result three
homogeneous linear equations in $\al_n,\al_{n-3}$ of which
the one for $m=2$ is trivial and
of which the other two yield $\al_n=\al_{n-3}=0$ unless $c=0,1,2$. Again,
together with \eqref{K9} and $f_{r,r,1}(c,d)\ne0$ this implies that all $\al_k$ are
zero if not $c=0,1,2$. For $c=0$ and for $c=2$ implying $d=1$ we already saw that
the solution space has dimension~1.

So the only remaining case to be considered is $(c,d)=(1,0)$.
We will show that
$q_n(x,y):=q_n^{(1,0)}(x,y)$, given by \eqref{K16},
clearly a polynomial of degree $n$, yields a (nonzero) element $q_n$ of
$\PP_{n-1}^{\perp}\left(  T_{1}\right)
\cap\PP_{n-1}^{\perp}\left(  K_{c,d}\right)$. By \eqref{K11} and
\eqref{K2} we see that
$(x,y)\mapsto P_{n+1}^{(1,1)}(1-2x)$
is an orthogonal polynomial of degree $n+1$ on $T_{1}$
with respect to the weight function
$x$ and that $(x,y)\mapsto P_{n+1}^{(1,1)}(1-2(x+y))=(-1)^{n}%
P_{n+1}^{(1,1)}(2(x+y)-1)$ is an orthogonal polynomial of degree $n+1$ on $T_{1}$
with respect to the weight function $1-x-y$. Then it holds
for any polynomial $r(x,y)$ of degree $<n$ that
\begin{multline*}
\int\!\!\!\int_{T_{1}}q_n(x,y)\,r(x,y)\,xy(1-x-y)\,dx\,dy\\
=\int\!\!\!\int_{T_{1}}P_{n+1}^{(1,1)}(1-2(x+y))\,(1-x-y)r(x,y)\,x\,dx\,dy\\
-\int\!\!\!\int_{T_{1}}P_{n+1}^{(1,1)}(1-2x)\,xr(x,y)\,(1-x-y)\,dx\,dy=0-0=0.
\end{multline*}
Hence $q_n\in\PP_{n}^{\perp}\left(  T_{1}\right)  $. Since $q_n$ is
invariant under the transformation $(x,y)\mapsto(x+y,-y)$, we have
$q_n\in\PP_{n}^{\perp}\left(  K_{1,0}\right) 
\cap\PP_{n}^{\perp}\left(  T_{1}\right) $. Since we already proved that in
this case the intersection has dimension at most one, we are finished.\quad
\endproof
\section{The intersection of $n$-th degree orthogonal polynomial spaces
for a triangle patch (case $n>1$)}
\label{Secn>1}
\begin{figure}[ptb]
\begin{center}
\includegraphics{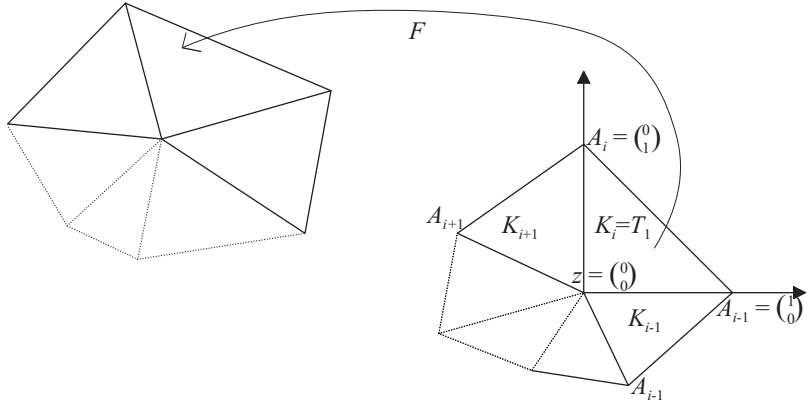}
\caption{Reference configuration. $K_{i}=T_{1}$ is the unit triangle,
$K_{i+1}$ is in the left half plane and $K_{i-1}$ in the lower half plane.}
\label{Fig2}
\end{center}
\end{figure}

Since polynomial spaces $\PP_{n}\left(D\right)$ are invariant under
affine coordinate transformations,
it suffices to prove Theorem \ref{MainConj}
for a
\emph{reference configuration}
with centre $z:=(0,0)$ and with one of the triangles, say $K_i$,
equal to the unit triangle $T_{1}$ given by \eqref{K12}.
Hence, the adjacent triangle $K_{i+1}$ to the left of $K_{i}$ lies in the left
half plane while the other one, i.e., $K_{i-1}$ lies in the lower half plane.
See Figure~\ref{Fig2}.

\begin{figure}[ptb]
\begin{center}
\includegraphics{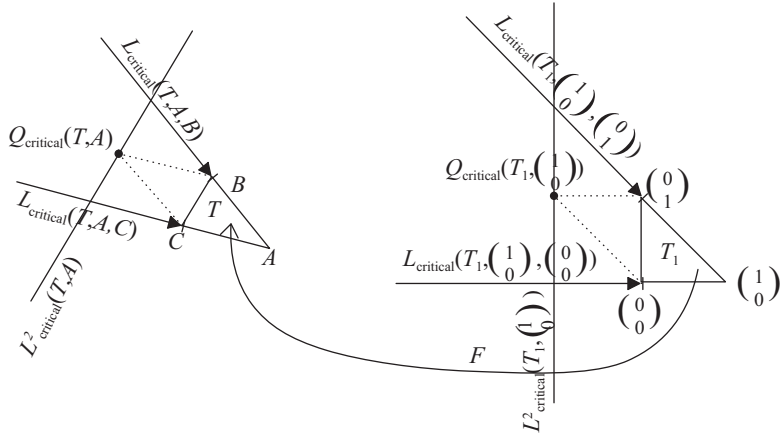}
\caption{Critical sets for the triangle $T$ and critical sets for the unit
triangle. }
\label{Fig3}
\end{center}
\end{figure}

In this Section we will prove the intersection property
Theorem \ref{MainConj}(b) for
polynomial degrees $n>1$.
First we will describe the exceptional cases in Theorem \ref{Theo2Triangles}
in terms of geometric quantities. For this we introduce the
``critical sets'' for a triangle; for an illustration see
Figure \ref{Fig3}.

\begin{definition}
For a triangle $T$ with vertices\footnote{As a convention we list the
vertices $A,B,C$ of a
triangle $T=\operatorname*{conv^o}\left(  A,B,C\right)$ always in the
counterclockwise ordering.} $A,B,C$, 
the \emph{critical sets} with respect to two vertices
$A,B$ are
\begin{align*}
&\Lcrit\left(  T,A,B\right):=\left\{
A+t\left(  B-A\right)  :t\geq1\right\}  ,\\
&Q_{\operatorname{critical}}\left(  T,A\right):=B+C-A,\\
&\Lcrit^2\left(  T,A\right):=\left\{
2B-A+t\left(C-B\right)  :t\in\RR\right\},\\
&\Lcrit^{n,\operatorname{tot}}\left(T,A\right)
:=\Lcrit\left(  T,A,B\right)  \cup
\Lcrit\left(  T,A,C\right)  \cup\left\{
Q_{\operatorname{critical}}\left(  T,A\right)  \right\}\;(n>2),\\
&\Lcrit^{2,\operatorname{tot}}\left(T,A\right)
:=\Lcrit\left(  T,A,B\right)  \cup
\Lcrit\left(  T,A,C\right)  \cup
\Lcrit^2\left(  T,A\right)\;(n=2).
\end{align*}
\end{definition}

Note that an orientation preserving affine map sending a triangle $T$ to a triangle
$\tilde T$ maps the critical sets of $T$ to the corresponding critical sets of $\tilde T$.
Also note that, for $n=2$,
$\Lcrit^2\left(T,A\right)$ contains
$Q_{\operatorname{critical}}\left(  T,A\right)$ and intersects with
$\Lcrit\left(T,A,B\right)$ and with
$\Lcrit\left(T,A,C\right)$.

\begin{proposition}
\label{LemLc1}Let $K_{1}=\operatorname*{conv^o}\left(  A,B,C\right)$
and $K_{2}=\operatorname*{conv^o}\left(  B,A,D\right) $ be two disjoint
triangles with common edge $AB$. Then for $n>1$
\[
\dim\left(\PP_{n-1}^{\perp}\left(  K_{1}\right)\cap
\PP_{n-1}^{\perp}\left(K_{2}\right) \right)  =\left\{
\begin{array}
[c]{cc}%
0 & \text{if }D\notin\Lcrit^{n,\operatorname{tot}%
}\left(  K_{1},C\right)  ,\\[\smallskipamount]
1 & \text{if }D\in\Lcrit^{n,\operatorname{tot}%
}\left(  K_{1},C\right)  .
\end{array}
\right.
\]
\end{proposition}%

\proof
We consider first the case that $K_{1}=T_{1}$, and apply Theorem \ref{Theo2Triangles}. Let $D=\left(\tfrac{-c}{d-c},\tfrac{1}{d-c}\right)$
(cf.~\eqref{defKparab}). The exceptional cases are given by

\begin{enumerate}
\item
$c=0$. Since $K_{1}$ and $K_{2}$ have disjoint interior this is
equivalent to
\[
D\in\Lcrit\big(T_{1},(0,1),(0,0)\big).
\]
\item
$d=1$. Again, taking into account that $K_{1}$ and $K_{2}$ have empty
open intersection we get that this case is equivalent to
\[
D\in\Lcrit\big(T_{1},(0,1),(1,0)\big).
\]
\item
$d-c=1$. This case contradicts the condition $K_{1}\cap K_{2}=\emptyset$
and, hence, cannot arise.
\item
$c=1$, $d=0$. This case is equivalent to
$D=Q_{\operatorname{critical}}\big(T_1,(0,1)\big) =\left(  1,-1\right)$.
\item
$n=2$ and $d=c-1$. Then $D=(s,-1)$ ($s\in\RR$), so this case is equivalent to
$D\in \Lcrit^2\big(T_1,(0,1)\big)$.
\end{enumerate}

The general case follows by employing an affine pullback of a general triangle
$K_{1}=\operatorname*{conv^o}\left(  A,B,C\right)  $ to $T_{1}$
such that $C$ is sent to $(0,1)$.\quad
\endproof
\bLP\goodbreak\noindent
{\bf Proof of Theorem \ref{MainConj}(b) for $n>1$}
\sLP
We use the numbering of triangles, edges, vertices in $\FST _{z}$ as in
Figure \ref{Fig1}. If there exists an edge $E_{i}$ with adjacent triangles
$K_{i}$, $K_{i+1}$ such that%
\[
A_{i+1}\notin\Lcrit^{n,\operatorname{tot}}\left(
K_{i},A_{i-1}\right)  \quad\vee\quad A_{i-1}\notin\FSL%
_{\operatorname{critical}}^{n,\operatorname{tot}}\left(  K_{i+1},A_{i+1}\right)
\]
we conclude from Proposition \ref{LemLc1} that%
\[
\PP_{n-1}^{\perp}\left(  K_{i}\right) 
\cap\PP_{n-1}^{\perp}\left(K_{i+1}\right)  =\left\{  0\right\}
\]
and the statement follows.

\begin{figure}[ptb]
\begin{center}
\includegraphics{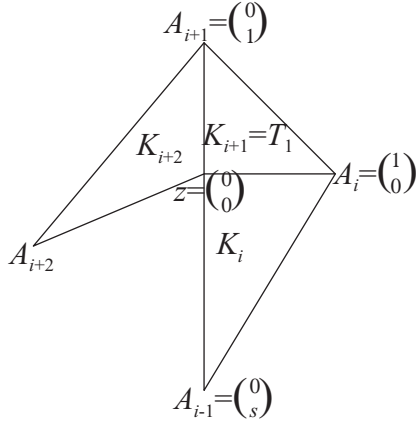}
\caption{Illustration of Case a) in the proof of Theorem \ref{MainConj}(b) for $n>1$.}%
\label{Figaligned}
\end{center}
\end{figure}

Hence, for the rest of the proof we always assume that%
\begin{equation}
\forall i\in\{1,\ldots,q\}\;
A_{i+1}\in\Lcrit^{n,\operatorname{tot}}\left(  K_{i},A_{i-1}\right)
\;\wedge\;
A_{i-1}\in\Lcrit^{n,\operatorname{tot}}
\left(  K_{i+1},A_{i+1}\right)  . \label{proofalwcrit}
\end{equation}
Clearly we can pick a vertex $A_i$ such that the inner angle at $A_i$,
(i.e., the angle
$\angle A_{i-1}A_iA_{i+1}$, seen from $z$) is less than $\pi$.
This property excludes that
$A_{i+1}\in
\Lcrit\left(K_{i},A_{i-1},A_i\right)$,
or equivalently
$A_{i-1}\in
\Lcrit\left(K_{i+1},A_{i+1},A_i\right)$.
The property is also preserved under affine maps.
We distinguish between the following cases.
\mLP
\textbf{Case a)}
$A_{i+1}\in
\Lcrit\left(K_{i},A_{i-1},z\right)$.\\
Without loss of generality we can work in
the reference situation (cf.~\mbox{Figure}~\ref{Figaligned})
that $K_{i+1}=T_{1}$, i.e., $A_{i+1}=\left(0,1\right)$,
$z=\left(0,0\right)$ and $A_{i-1}=\left(0,s\right)$ for some $s<0$. From Theorem \ref{Theo2Triangles}, in particular from \eqref{K18},
it follows that
$\PP_{n-1}^{\perp}\left(K_{i+1}\right)\cap
\PP_{n-1}^{\perp}\left(K_{i}\right)$
is spanned by the polynomial
$q(x,y):=P_{n}^{(1,3)}\left(1-2x\right)$.
Note that the adjacent triangle $K_{i+2}$ left to $T_{1}$ lies in the left
half plane. Hence $q$ is either
positive on $K_{i+2}$ or negative, by which it cannot be orthogonal to all
constant functions on $K_{i+2}$. We conclude that
\[
\PP_{n-1}^{\perp}\left(  K_{i+2}\right)  \cap\PP_{n-1}^{\perp
}\left(  K_{i+1}\right)  \cap\PP_{n-1}^{\perp}\left(  K_{i}\right)
=\left\{  0\right\}  .
\]
\textbf{Case b)}\quad
$A_{i+1}=Q_{\operatorname{critical}}\left(  K_{i},A_{i-1}\right)$
and $A_{i-2}\in\Lcrit\left(
K_{i},A_{i},z\right)$.\\
Then
$A_i\in
\Lcrit\left(K_{i-1},A_{i-2},z\right)$.
So as in case a), but now with $K_i=T_1$, we conclude that
\[
\PP_{n-1}^{\perp}\left(  K_{i-1}\right)\cap
\PP_{n-1}^{\perp}\left(  K_{i}\right)\cap
\PP_{n-1}^{\perp}\left(  K_{i+1}\right)
=\left\{  0\right\}  .
\]
\textbf{Case c)}\quad
$A_{i+1}=Q_{\operatorname{critical}}\left(  K_{i},A_{i-1}\right)$ and
$A_{i-2}\in\Lcrit\left(  K_{i} ,A_{i},A_{i-1}\right)$.
\begin{figure}[ptb]
\begin{center}
\includegraphics{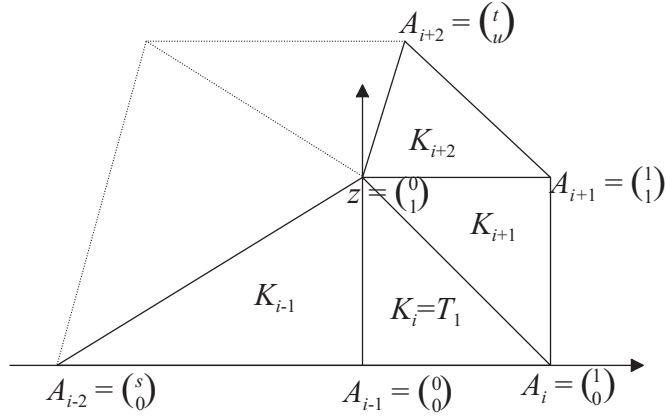}
\caption{Illustration of the geometric argument for case c). The points
$A_{i-2},A_{i-1},A_{i}$ are collinear and $\overline{K_{i}\cup K_{i+1}}$ form
a parallelogram.}
\label{Fig6}
\end{center}
\end{figure}
\\
Without loss of generality we can work in the
reference configuration (see Figure \ref{Fig6}) that
$K_{i}=T_{1}$ with $A_{i-1}=\left(0,0\right)$,
$A_{i}=\left(1,0\right)$,
$z=\left(0,1\right)$. Then
$K_{i+1}=\operatorname*{conv^o}\left((1,0),(1,1),(0,1)\right)$,
$K_{i-1}=\operatorname*{conv^o}\left((0,0),(0,1),(s,0)\right)$
for some $s<0$, and
$K_{i+2}=\operatorname*{conv^o}\left((0,1),(1,1),(t,u)\right)$
with $u>1$.
From Theorem \ref{Theo2Triangles},
in particular from \eqref{K20}, we conclude that
$\PP_{n-1}^{\perp}\left(K_{i}\right)\cap
\PP_{n-1}^{\perp}\left(K_{i-1}\right)$ is spanned by
the polynomial  $r(x,y):=P_{n}^{\left(1,3\right)}\left(1-2y\right)$.
By arguing as in Case a) we conclude that $r$ is not changing sign in
$K_{i+2}$ (since $y\geq1$).
Hence
\[
\PP_{n-1}^{\perp}\left(  K_{i-1}\right)  \cap\PP_{n-1}^{\perp}
\left(  K_{i}\right)  \cap\PP_{n-1}^{\perp}\left(  K_{i+1}\right)
\cap\PP_{n-1}^{\perp}\left(  K_{i+2}\right)  =\left\{  0\right\}  .
\]
\textbf{Case d)}\quad$A_{i+1}=Q_{\operatorname{critical}}\left(  K_{i},A_{i-1}\right)$
and
$A_{i-2}=Q_{\operatorname{critical}}\left(K_{i},A_i\right)$.\\
Without loss of generality we can consider the reference situation
that $K_i=T_1$ with $z=(0,0)$, $A_{i-2}=(1,-1)$, $A_{i-1}=(1,0)$,
$A_i=(0,1)$, $A_{i+1}=(-1,1)$.
From Theorem \ref{Theo2Triangles}
we conclude that
$\PP_{n-1}^{\perp}\left(K_{i}\right)\cap
\PP_{n-1}^{\perp}\left(K_{i-1}\right)$ is spanned by
the polynomial $q_n(x,y):=q_n^{(1,0)}(x,y)$ given by \eqref{K16}, and that
$\PP_{n-1}^{\perp}\left(K_{i+1}\right)\cap
\PP_{n-1}^{\perp}\left(K_i\right)$ is spanned by
the polynomial $q_n(y,x)$.
Since these two polynomials are linearly independent (compare the highest degree
part of both polynomials), we have shown that
\[
\PP_{n-1}^{\perp}\left(  K_{i-1}\right)\cap
\PP_{n-1}^{\perp}\left(  K_{i}\right)\cap
\PP_{n-1}^{\perp}\left(  K_{i+1}\right)
=\left\{  0\right\}.
\]
\textbf{Case e)}\quad
$A_{i+1}\in
\Lcrit^2\left(K_i,A_{i-1}\right)$
and $A_{i-2}\in
\Lcrit\left(  K_{i} ,A_{i},A_{i-1}\right)$.\\
We can assume that $K_i=T_1$ with $z=(0,0)$. Then $A_{i+1}=(-1,c)$.
From \eqref{K21} we have that
$\PP_{n-1}^{\perp}\left(K_{i+1}\right)\cap
\PP_{n-1}^{\perp}\left(K_i\right)$ is spanned by
the polynomial $q_2^{(c,c-1)}(y,x)$ and from \eqref{K19} that
$\PP_{n-1}^{\perp}\left(K_{i}\right)\cap
\PP_{n-1}^{\perp}\left(K_{i-1}\right)$ is spanned by
the polynomial $P_n^{(1,3)}(2(x+y)-1)$. A computation shows that
these two polynomials are linearly dependent iff $c=2$. But then the
inner angle at $A_i$ equals $\pi$, which we excluded.
\\[\smallskipamount]
\textbf{Case f)}\quad
$A_{i+1}\in
\Lcrit^2\left(K_i,A_{i-1}\right)$
and $A_{i-2}\in
\Lcrit\left(K_{i},A_{i},z\right)$.\\
Again assume that $K_i=T_1$ with $z=(0,0)$. Then $A_{i-2}=(0,s)$ for some $s<0$.
We can argue as in case a), but now with $K_i=T_1$.
\\[\smallskipamount]
\textbf{Case g)}\quad
$A_{i+1}\in
\Lcrit^2\left(K_i,A_{i-1}\right)$
and $A_{i-2}\in
\Lcrit^2\left(K_i,A_i\right)$.\\
Again assume that $K_i=T_1$ with $z=(0,0)$.
Then $A_{i+1}=(-1,c)$ and $A_{i-2}=(d,-1)$ for some $c,d\in\RR$.
From \eqref{K21} we have that
$\PP_{n-1}^{\perp}\left(K_{i+1}\right)\cap
\PP_{n-1}^{\perp}\left(K_i\right)$ is spanned by
the polynomial $q_2^{(c,c-1)}(y,x)$ and that
$\PP_{n-1}^{\perp}\left(K_{i}\right)\cap
\PP_{n-1}^{\perp}\left(K_{i-1}\right)$ is spanned by
the polynomial $q_2^{(d,d-1)}(x,y)$. A computation shows that these
two polynomials are linearly dependent iff $c=d=-3$ or $c=d=2$.
But in the first case the triangles $K_{i-1}$ and
$K_{i+1}$ intersect, and in the second case the inner angle at $A_i$
equals $\pi$.\quad\endproof

\section{The intersection of $n$-th degree orthogonal polynomial spaces
for a triangle patch (case $n=1$)}
\label{Secn=1}
In this Section we will prove the intersection property
Theorem \ref{MainConj}(b)
polynomial degrees $n=1$.
We need two simple lemmas.

\begin{lemma}
\label{K23}
Let $\FST :=\left\{  K_{i}:1\leq i\leq q\right\}  $ denote a triangle
patch around $z\in\mathbb{R}^{2}$.
Then there are $K_{i-1}$, $K_i$, $K_{i+1}$ of which the barycenters are not
collinear.
\end{lemma}

\proof
Let $M_i$ be the barycenter of $K_i$. We will show that the points $M_i$
($i=1,\ldots,q$) cannot be collinear. This will imply the statement of the Lemma.

\begin{figure}[ptb]
\begin{center}
\includegraphics{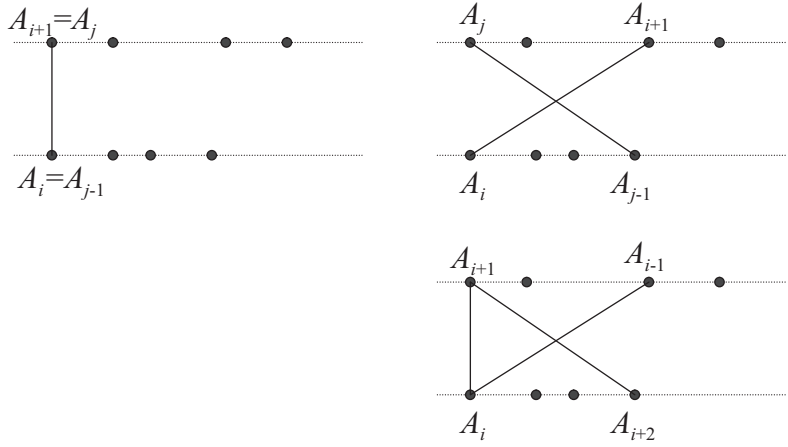}
\caption{Illustration of the last part of the proof of Lemma \ref{K23}  ($q$ even).
$A_{i+1}=A_{j}$ (left picture)  because otherwise the top right picture gives
a contradiction. 
But now the bottom
picture gives a contradiction.}
\label{FigAlignedVertices}
\end{center}
\end{figure}

We may choose $z=(0,0)$. Suppose that the points $M_i=\tfrac13(A_{i-1}+A_i)$
($i=1,\ldots,q$) are collinear. Then all vectors
$3(M_i-M_{i-1})=A_i-A_{i-2}$ are proportional. If $q$ is odd then this implies
that all vertices $A_i$ are collinear, which is impossible.
If $q$ is even then the set of vertices $A_1,A_3,\ldots,A_{q-1}$ 
and the set of vertices $A_2,A_4,\ldots,A_q$ are both collinear and the two
collinear sets lie on parallel lines. Since all vertices cannot be
collinear, these two parallel lines have to be distinct.
After applying an affine linear map we may assume that one of the lines is
$y=0$ with $A_i=(0,0)$ for some $i$ and with all other vertices on this line having
coordinates $(x,0)$ with $x>0$, and that the other line is
$y=1$ with $A_j=(0,1)$ for some $j\ne i$ and with all other vertices on this line having
coordinates $(x,1)$ with $x>0$.
First we show that $j=i+1$. Indeed, if
$j\ne i+1$ then
the edge connecting $A_j$ and $A_{j-1}$
will cross the edge connecting $A_i$ and $A_{i+1}$, which is not allowed.
Thus $A_{i+1}=(0,1)$.
But now the edge connecting $A_{i+1}$ and $A_{i+2}$
will cross the edge connecting $A_i$ and $A_{i-1}$,
which is not allowed (see Figure \ref{FigAlignedVertices} for this
last part of the proof, where we successively arrive twice at a contradiction).
Thus we cannot have two collinear sets of vertices if $q$ is even.\quad\endproof
\begin{lemma}
\label{K24}
Let $K$ be a triangle with barycenter $M$. Let $\mu$ be a finite measure on $K$
which is invariant under all affine transformations mapping $K$ onto itself
(these form a group ismorphic with $S_3$). Then
\[
\int_K p\,d\mu=p(M)\,\mu(K)\quad\mbox{for all $p\in\PP_1(K)$.}
\]
This holds in particular if $d\mu(x,y)=\om(x,y)\,dx\,dy$ with $\om$ the product
of the barycentric coordinates for $K$.
\end{lemma}

\proof
Let $A,B,C$ be the vertices of $K$.
Since the assertion is trivial for constant functions,
it is sufficient to prove the property for affine linear functions $p$
vanishing on one of the medians $AM$, $BM$, $CM$. Suppose $p$ vanishes
on $AM$.  Then the function $p-p(M)$ is sent to its opposite
under the affine map fixing $A$ and interchanging $B$ and $C$
(check this for the reference triangle $T$
with $A=(0,0)$, $B=(1,0)$, $C=(0,1)$).
Hence $\int_K (p-p(M))\,d\mu=0$. \quad\endproof
\bLP
{\bf Proof of Theorem \ref{MainConj}(b) for $n=1$}
\sLP
Let $\FST :=\left\{  K_{i}:1\leq i\leq q\right\}$ denote a triangle
patch around a point $z\in\mathbb{R}^{2}$
and let $\Om:=\cup_{i=1}^{q}\overline{K_{i}}$.
By Lemma \ref{K23} there are $K_{i-1},K_i,K_{i+1}$ such that their
barycenters $M_{i-1},M_i,M_{i+1}$ are not collinear.
Now suppose that $u\in\PP_{1}(\Om)$ and
$u\in\PP_0^\perp(K_j)$ for $j=i-1,i,i+1$.
By Lemma \ref{K24} we have for $j\in\{i-1,i,i+1\}$ that
\[
0=\frac{\int_{K_j}u(x,y)\,\om_j(x,y)\,dx\,dy}
{\int_{K_j}\om_j(x,y)\,dx\,dy}=u(M_j).
\]
Since the affine linear function $u$ vanishes on three points which
are not collinear, $u$ is identically zero.\quad\endproof

\begin{remark}
If we take the triangle patch such that
$K_i=T_1$ with $z=(0,0)$, $A_{i-1}=(1,0)$, $A_i=(0,1)$,
$A_{i-2}=(\tfrac{-c}{d-c},\tfrac1{d-c})$ and
$A_{i+1}=(\tfrac1{d'-c'},\tfrac{-c'}{d'-c'})$,
then we see from \eqref{K22} that
\[
\PP_{0}^{\perp}\left(  K_{i-1}\right)\cap
\PP_{0}^{\perp}\left(  K_{i}\right)\cap
\PP_{0}^{\perp}\left(  K_{i+1}\right)
\ne\left\{  0\right\}  .
\]
implies that the polynomials $q_1^{(c,d)}(x,y)$ and $q_1^{(c',d')}(y,x)$ are
multiples of each other.
A computation shows that then
\[
c'=k(c-d+1),\quad
d'=k(1-d)+1\qquad(0\ne k\in\RR).
\]
Hence $A_{i+1}=(\tfrac1{1-kc},\tfrac{k(d-c-1)}{1-kc})$.
Then $A_{i+1}-A_{i-1}$ and $A_i-A_{i-2}$ are proportional.
By the Proof of Lemma \ref{K23} this gives the collinearity of the barycenters
of $K_{i-1}$, $K_i$ and $K_{i+1}$. Thus we have shown once more that,
if the barycenters of these three triangles are not collinear, then
$\PP_{0}^{\perp}\left(  K_{i-1}\right)\cap
\PP_{0}^{\perp}\left(  K_{i}\right)\cap
\PP_{0}^{\perp}\left(  K_{i+1}\right)
=\left\{  0\right\}$.
\end{remark}

\section{Injectivity for the polynomial projection\\
operator: some follow-up}
\label{K33}
Recall Theorem \ref{MainConj}(c) about the injectivity of
the polynomial projection operator 
$\Pi_n^{\FST}\colon\PP_{n}(\Om)\to\PP_{n-1}(\FST)$,
Let $\left\Vert
\cdot\right\Vert _{\FST }:=\left(  \cdot,\cdot\right)  _{\FST %
}^{1/2}$ and define
\begin{align}
c_n'(\FST)&:=\inf_{v\in\PP_n(\Om)\backslash\{0\}}
\frac{(v,\Pi_n^\FST v)_\FST}{(v,v)_\FST}
=\inf_{v\in\PP_n(\Om)\backslash\{0\}}
\frac{\|\Pi_n^\FST v\|_\FST^2}{\|v\|_\FST^2}>0,
\label{K26}\\
c_n''(\FST)&:=\inf_{v\in\PP_n(\Om)\backslash\{0\}}
\frac{\|v\|_\FST^2}{\|v\|_{L^2(\Om)}^2}>0,
\label{K27}\\
\check c_n(\FST)&:=\inf_{v\in\PP_n(\Om)\backslash\{0\}}
\frac{(v,\Pi_n^\FST v)_\FST}{(v,v)_{L^2(\Om)}}\ge
c_n''(\FST)c_n'(\FST)>0.
\label{K28}
\end{align}
The inequality in \eqref{K26} follows from the injectivity of $\Pi_n^\FST$,
while the inequality in \eqref{K27} is a consequence of the equivalence
of all norms on a finite dimensional space.

\begin{remark}
We can expand the squared norms in \eqref{K26} and \eqref{K27}.
Let $\La_i$ be the orientation preserving affine linear map
of $T_1$ onto $K_i$ which sends $(0,0)$ to $z$.
Let $\{p_m\}_{m=1,\ldots,n(n+1)/2}$ be an orthonormal basis for $\PP_{n-1}(T_1)$
with respect to the weight function $\om$ on $T_1$. Then
\begin{align*}
\big\|\Pi_n^\FST v\big\|_\FST^2=
\sum_{i=1}^q \det\La_i \sum_{m=1}^{\half n(n+1)}
&\left(\int_{T_1} v\circ\La_i\;p_m\,\om\right)^2,\\
\|v\|_\FST^2=\sum_{i=1}^q \det\La_i \int_{T_1}(v\circ\La_i)^2\,\om,\qquad&
\|v\|_{L^2(\Om)}^2=\sum_{i=1}^q \det\La_i \int_{T_1}(v\circ\La_i)^2.
\end{align*}
Hence $c_n''=c_n''(\FST)$ is independent of the choice of the patch $\FST$
and it equals
\begin{equation}
c_n''=\inf_{v\in\PP_n(T_1)\backslash\{0\}}
\frac{\int_{T_1} v(x,y)^2\,x\,y\,(1-x-y)\,dx\,dy}
{\int_{T_1} v(x,y)^2\,dx\,dy}\,.
\label{K37}
\end{equation}
\end{remark}

\begin{remark}
From \cite[Prop.~3.46]{Verfuerth} there follows an $n$-explicit lower bound
for $c_n''$ as given by \eqref{K37}:
$c_n''\ge C/(n+1)^4$
with a fixed $n$-independent constant $C>0$.
\end{remark}

\begin{figure}[ptb]
\begin{center}
\includegraphics{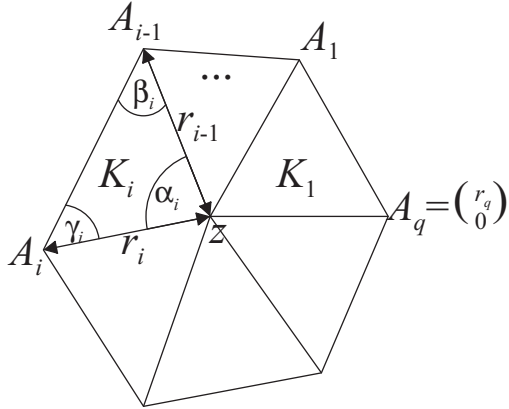}
\caption{Illustration of angles $\al_1,\be_i,\ga_i$, and radii
$r_i$ for a triangle patch.}
\label{Fig_patch_para}
\end{center}
\end{figure}

Since the quotients of integrals in \eqref{K26}--\eqref{K28} are invariant
under affine linear maps, we might restrict to the case that $K_1$ is the
reference triangle $T_1$. But in view of the numerical applications, we consider
only translations and rotations of triangle patches. Thus we
restrict in the following to the case that $z=0$ and $A_q=(r_q,0)$ for some $r_q>0$.
The further data determining $\FST$ are, for $i=1,\ldots,q$,
the angles $\al_i,\be_i,\ga_i$ of the triangle $K_i$
at $z,A_{i-1},A_i$, respectively, together with the length $r_i$ of the edge connecting
$z$ with $A_i$
(see Figure \ref{Fig_patch_para}).
Evident constraints on these numbers are that
$\al_i+\be_i+\ga_i=\pi$ and $\al_1+\cdots+\al_q=2\pi$. But $\FST$ would
already be completely determined by $\al_2,\ldots,\al_q$, $\be_2,\ldots,\be_q$
and $r_q$, or by $\al_2,\ldots,\al_q$ and $r_1,\ldots,r_q$.
The map
$(\al_2,\ldots,\al_q,\be_2,\ldots,\be_q,r_q)\leftrightarrow
(\al_2,\ldots,\al_q,r_1,\ldots,r_q)$
is continuous in both directions, as can be seen from the following identities
obtained by a combination of the sine rule and the cosine rule for the triangle $K_i$\,:
\begin{equation}
\frac{\sin\al_i}{\sqrt{r_{i-1}^2+r_i^2-2r_{i-1}r_i\cos\al_i}}
=\frac{\sin\be_i}{r_i}=\frac{\sin(\pi-\al_i-\be_i)}{r_{i-1}}\,.
\label{K34}
\end{equation}

The quotients of integrals in \eqref{K26}--\eqref{K28} will depend continuously
on $v$ and the data of $\FST$. Therefore, the three constants in
\eqref{K26}--\eqref{K28} will depend continuously on the data of $\FST$
and they will remain bounded away from zero if we let range the data of $\FST$
over a compact set. To fix a compact set, choose $\de\in(0,\pi/3]$ and
$\rho>0$.

\begin{definition}\label{K29}
The compact set of triangle patches $X_{q,\de,\rho}$ consists of all $\FST$
with $\al_i,\be_i,\ga_i\ge\de$ and $r_i\ge\rho$ ($i=1,\ldots,q$).
Furthermore define
\begin{equation}
\tilde c_n(q,\de,\rho):=\inf_{\FST\in X_{q,\de,\rho}} \check c_n(\FST).
\label{K36}
\end{equation}
\end{definition}

By the second equality in \eqref{K34} we see that, for given $\de$ and $q$ there
exists $C>0$ such that $\rho\le r_i\le C\rho$ ($i=1,\ldots,q)$ if
$\FST\in X_{q,\de,\rho}$.
Since the quotients of integrals in \eqref{K26}--\eqref{K28} are invariant under
dilations, $\tilde c_n(q,\de,\rho)$ will be independent of $\rho$.
Since necessarily $q\de\le 2\pi$, $X_{q,\de,\rho}$ is non-empty for only finitely
many values of $q$. We conclude:

\begin{theorem}
$\displaystyle \inf_{\rho>0,\;q\ge3}\tilde c_n(q,\de,\rho)>0$.
\label{K35}
\end{theorem}

\bibliographystyle{abbrv}
\bibliography{poly_intersect}

\end{document}